\documentclass[12pt, a4paper]{article}
\usepackage[latin1]{inputenc}
\usepackage[english]{babel}
\usepackage{amsmath}
\usepackage[english]{babel}
\usepackage{indentfirst}
\usepackage{amstext}
\usepackage{enumerate}
\usepackage{amsfonts}
\usepackage{textcomp}
\usepackage{amssymb}
\usepackage[dvips]{graphicx}
\usepackage{setspace}
\usepackage[all]{xy}
\usepackage{color}

\newcommand {\demo}{\hskip -0.6cm{\bf Proof.  }}
\newcommand {\fim}{\hfill{$\square$}\vskip 1pc}

\newcommand {\N}{\mathbb{N}}

\newcommand {\F}{\mathbb{F}}

\newcommand {\GG}{\mathcal{G}}

\newcommand{\id}{\mathrm{id}}

\newtheorem{teorema}{Theorem}[section]
\newtheorem{lema}[teorema]{Lemma}
\newtheorem{corolario}[teorema]{Corollary}
\newtheorem{definicao}[teorema]{Definition}
\newtheorem{proposicao}[teorema]{Proposition}

\newtheorem{remark}[teorema]{Remark}

\begin{document}
\onehalfspace

\title{Representations and the reduction theorem for ultragraph Leavitt path algebras}

\author{Daniel Gon\c{c}alves\footnote{This author is partially supported by Conselho Nacional de Desenvolvimento Cient\'{i}fico e Tecnol\'{o}gico - CNPq.}  \ and Danilo Royer}
\maketitle

\begin{abstract}
In this paper we study representations of ultragraph Leavitt path algebras via branching systems and, using partial skew ring theory, prove the reduction theorem for these algebras. We apply the reduction theorem to show that ultragraph Leavitt path algebras are semiprime and to completely describe faithfulness of the representations arising from branching systems, in terms of the dynamics of the branching systems. Furthermore, we study permutative representations and provide a sufficient criteria for a permutative representation of an ultragraph Leavitt path algebra to be equivalent to a representation arising from a branching system. We apply this criteria to describe a class of ultragraphs for which every representation (satisfying a mild condition) is permutative and has a restriction that is equivalent to a representation arising from a branching system. 

\end{abstract}

\let\thefootnote\relax

\section{Introduction}

The study of algebras associated to combinatorial objects is a mainstream area in Mathematics, with connections with symbolic dynamics, wavelet theory and graph theory, to name a few. Among the most studied algebras arising from combinatorial objects we find Cuntz-Krieger algebras, graph C*-algebras, Leavitt path algebras (algebraic version of graph C*-algebras) and algebras associated to infinite matrices (the so called Exel-Laca algebras). Aiming at an unified approach to graph C*-algebras and Exel-Laca algebras Mark Tomforde introduced ultragraphs in \cite{Tomforde:JOT03}. One of the advantages of dealing with ultragraphs arises from the combinatorial picture, very similar to graphs, available to study Exel-Laca algebras. Also, new examples appear, as the class of ultragraph algebras is strictly larger then the class of graph and Exel-Laca algebras (both in the C*-context and in the purely algebraic context), although in the C*-algebraic context these three classes agree up to Morita equivalence, see \cite{KMST}.

Ultragraphs can be seen as graphs for which the range map takes values over the power set of the vertices, that is, the range is a subset of the set of vertices. Ultragraphs can also be seen as labelled graphs where it is only possible to label edges with the same source. This restriction makes it much simpler to deal with objects (such as algebras and topological spaces) associated to ultragraphs than with objects associated to labelled graphs, but still interesting properties of labelled graphs present themselves in ultragraphs. This is the case, for example, in the study of Li-Yorke chaos for shift spaces over infinite alphabets (see \cite{LiYorke}). In fact, ultragraphs are key in the study of shift spaces over infinite alphabets, see \cite{GRultrapartial, GRultra, GSCSC}. Furthermore, the KMS states associated to ultragraph C*-algebras are studied in \cite{GD} and the connection of ultragraph C*-algebras with the Perron-Frobenius operator is described in \cite{ultragraphrep}, where the theory of representations of ultragraph C*-algebras (arising from branching systems) is also developed. It is worth mentioning that many of the results in \cite{ultragraphrep} are lacking an algebraic counter-part and this is partially the goal of this paper.

Over the years, many researches dedicated efforts to obtain analogues of results in operator theory in the purely algebraic context, and to understand the relations between these results. For example, Leavitt path algebras, see \cite{AAS, AA, AMP}, were introduced as an algebraisation of graph C*-algebras and Cuntz-Krieger algebras. Later, Kumjian-Pask algebras, see \cite{PCHR}, arose as an algebraisation of higher rank graph C*-algebras. Partial skew group rings were studied as algebraisation of partial crossed products, see \cite{MishaExel} and  Steinberg algebras were introduced in \cite{BS, CFST} as an algebraisation of the groupoid C*-algebras first studied by Renault. Very recently the algebraisation of ultragraph C*-algebras, called ultragraph Leavit path algebras was defined, see \cite{leavittultragraph}. Similarly to the C*-algebraic setting, ultragraph Leavit path algebras generalize the Leavitt path algebras and the algebraic version of Exel-Laca algebras and provide for examples that are neither Leavitt path algebras nor Exel-Laca algebras. 

As we mentioned before our goal is to study ultragraph Leavitt path algebras. Our first main result is the reduction theorem, which for Leavitt path algebras was proved in \cite{Malaga}. This result is fundamental in Leavitt path algebra theory and it is also key in our study of representations of ultragraph Leavitt path algebras (we also use it to prove that ultragraph Leavitt path algebras are semiprime). The study of representations of algebras associated to combinatorial objects is a subject of much interest. For example, representations of Leavitt path algebras were studied in \cite{AR, Chen, HR, Ranga}, of 
Kumjian-Pask in \cite{PCHR}, of Steinberg algebras in \cite{AHLS, BCFS}. Representations of various algebras, in connection with branching systems, were studied in \cite{CG, CNPS, MR3404559, FGKP, GLR, GLR2, GR3, represleavitt, unitaryequivalence, MR2903145, GR4}. To describe the connections of representations of ultragraph Leavitt path algebras with branching systems is the second goal of this paper. In particular, we will give a description of faithful representations arising from branching systems, will define permutative representations and show conditions for equivalence between a given representation and representations arising from branching systems (these are algebraic versions of the results in \cite{ultragraphrep}).

The paper is organized as follows: The next section is a brief overview of the definitions of ultragraphs and the associated ultragraph Leavitt path algebra. In Section~3 we prove the reduction theorem, using partial skew group ring theory and the grading of ultragraph Leavitt path algebras by the free group on the edges (obtained from the partial skew ring characterization). We notice that the usual proof of the reduction theorem for Leavitt path algebras does not pass straighforwardly to ultragraph Leavitt path algebras and so we provide a proof using partial skew group ring theory (in the case of a graph our proof is an alternative proof of the reduction theorem for Leavitt path algebras). In Section~4 we define branching systems associated to ultragraphs and show how they induce a representation of the algebra. The study of faithful representations arising from branching systems is done in Section~5. Finally, in Section~6 we define permutative representations and study equivalence of representations of an ultragraph Leavitt path algebra with representations arising from branching systems. 

Before we proceed we remark that we make no assumption of countability on the ultragraphs, and hence the results we present generalize, to the uncountable graph case, the results for Leavitt path algebras presented in \cite{represleavitt}.

\section{Preliminaries}

We start this section with the definition of ultragraphs.

\begin{definicao}[{\cite[Definition~2.1]{Tomforde:JOT03}}]\label{def of ultragraph}
An \emph{ultragraph} is a quadruple $\mathcal{G}=(G^0, \mathcal{G}^1, r,s)$ consisting of two sets $G^0, \mathcal{G}^1$, a map $s:\mathcal{G}^1 \to G^0$ and a map $r:\mathcal{G}^1 \to P(G^0)\setminus \{\emptyset\}$, where $P(G^0)$ stands for the power set of $G^0$.
\end{definicao}

Before we define the algebra associated with an ultragraph, we need a notion of "generalized vertices". This is the content of the next definition.

\begin{definicao}[{\cite{Tomforde:JOT03}}]\label{defofultragraph}
Let $\mathcal{G}$ be an ultragraph. Define $\mathcal{G}^0$ to be the smallest subset of $P(G^0)$ that contains $\{v\}$ for all $v\in G^0$, contains $r(e)$ for all $e\in \mathcal{G}^1$, and is closed under finite unions and finite intersections.
\end{definicao}

Notice that since $\mathcal{G}^0$ is closed under finite intersections, the emptyset is in $\mathcal{G}^0$. We also have the following helpful description of the set of generalized vertices $\mathcal{G}^0$.

\begin{lema}[{\cite[Lemma~2.12]{Tomforde:JOT03}}]\label{concrete des of G^0}
Let $\mathcal{G}$ be an ultragraph. Then
\begin{align*}
\mathcal{G}^0=\Big\{\Big(\bigcap_{e\in X_1}r(e)\Big) \cup\dots\cup\Big(\bigcap_{e\in X_n}r(e)\Big) \cup F&:X_i \text{'s are finite subsets of } \mathcal{G}^1,
\\&F \text{ is a finite subset of } G^0\Big\}.
\end{align*}
\end{lema}

We can now define the ultragraph Leavitt path algebra associated to an ultragraph $\GG$. 

\begin{definicao}[\cite{leavittultragraph}]\label{def of ultragraph algebra}
Let $\mathcal{G}$ be an ultragraph and $R$ be a unital commutative ring. The Leavitt path algebra of $\mathcal{G}$, denoted by $L_R(\mathcal{G})$ is the universal $R$ with generators $\{s_e,s_e^*:e\in \mathcal{G}^1\}\cup\{p_A:A\in \mathcal{G}^0\}$ and relations
\begin{enumerate}
\item $p_\emptyset=0,  p_Ap_B=p_{A\cap B},  p_{A\cup B}=p_A+p_B-p_{A\cap B}$, for all $A,B\in \mathcal{G}^0$;
\item $p_{s(e)}s_e=s_ep_{r(e)}=s_e$ and $p_{r(e)}s_e^*=s_e^*p_{s(e)}=s_e^*$ for each $e\in \mathcal{G}^1$
\item $s_e^*s_f=\delta_{e,f}p_{r(e)}$ for all $e,f\in \mathcal{G}$
\item $p_v=\sum\limits_{s(e)=v}s_es_e^*$ whenever $0<\vert s^{-1}(v)\vert< \infty$.
\end{enumerate}
\end{definicao}

To prove the reduction theorem in the next section we need the characterization of Leavitt ultragraph algebras as partial skew rings. Therefore we recall this description below (as done in \cite{crossedproduct}).

\subsection{Leavitt ultragraph path algebras as partial skew rings}

We start setting up some notation. A finite path is either an element of $\mathcal{G}^0$ or a sequence of edges $e_1...e_n$, with length $|e_1...e_n|=n$, and such that $s(e_{i+1})\in r(e_i)$ for each $i\in \{0,...,n-1\}$. An infinite path is a sequence $e_1e_2e_3...$, with length $|e_1e_2...|=\infty$, such that $s(e_{i+1})\in r(e_i)$ for each $i\geq 0$. The set of finite paths in $\mathcal{G}$ is denoted by $\mathcal{G}^*$, and the set of infinite paths in $\mathcal{G}$ is denoted by $\mathfrak{p}^\infty$. We extend the source and range maps as follows: $r(\alpha)=r(\alpha_{|\alpha|})$, $s(\alpha)=s(\alpha_1)$ for $\alpha\in \mathcal{G}^*$ with $0<|\alpha|<\infty$, $s(\alpha)=s(\alpha_1)$ for each $\alpha\in \mathfrak{p}^\infty$, and $r(A)=A=s(A)$ for each $A\in \mathcal{G}^0$. An element $v\in G^0$ is a sink if $s^{-1}(v) = \emptyset$, and we denote the set of sinks in $G^0$ by $G^0_s$. We say that $A\in \mathcal{G}^0$ is a sink if each vertex in $A$ is a sink. 

Define the set $$X=\mathfrak{p}^\infty\cup\{(\alpha,v): \alpha \in \mathcal{G}^*, |\alpha|\geq 1, v \in G^0_s \cap r(\alpha) \}\cup\{(v,v): v\in G^0_s \}  .$$

We extend the range and source maps to elements $(\alpha,v)\in X$ by defining $r(\alpha,v)=v$ and $s(\alpha,v)=s(\alpha)$. Furthermore, we extend the length map to the elements $(\alpha,v)$ by defining $|(\alpha,v)|:=|\alpha|$.

The group acting on the space $X$ is the free group generated by $\mathcal{G}^1$, which we denote by $\F$. Let $W\subseteq \F$ be the set of paths in $\GG^*$ with strictly positive length.

Now we define the following sets:
\begin{itemize}
\item for $a\in W$, let $X_a=\{x\in X:x_1..x_{|a|}=a\}$;

\item for $b\in W$, let $X_{b^{-1}}=\{x\in X:s(x)\in r(b)\}$;
\item for $a,b\in W$ with $r(a)\cap r(b)\neq \emptyset$, let $$X_{ab^{-1}}=\left \{x\in X:|x|>|a|, \,\,\,  x_1...x_{|a|}=a \text{ and }s(x_{|a|+1})\in r(b)\cap r(a)\right\}\bigcup$$ $$\bigcup\left\{(a,v) \in X:v\in r(a)\cap r(b)\right\};$$
\item for the neutral element $0$ of $\F$, let $X_0=X$;
\item for all the other elements $c$ of $\F$, let $X_c=\emptyset$.
\end{itemize}

Define, for each $A\in \mathcal{G}^0$ and $b\in W$, the sets 
$$X_A=\{x\in X:s(x)\in A \}$$ and 
$$X_{bA}=\{x\in X_b:|x|>|b| \text{ and } s(x_{|b|+1})\in A\}\cup\{(b,v)\in X_b:v\in A\}.$$

We obtain a partial action of $\F$ on $X$ by defining the following bijective maps:
\begin{itemize}
\item for $a\in W$ define $\theta_a:X_{a^{-1}}\rightarrow X_a$ by $$\theta_a(x)=\left\{\begin{array}{ll}
ax & \text{ if } |x|=\infty,\\
(a\alpha,v) & \text{ if } x=(\alpha,v),\\
(a,v) & \text{ if } x=(v,v);

\end{array} \right .$$

\item for $a\in W$ define $\theta_a^{-1}:X_a\rightarrow X_{a^{-1}}$ as being the inverse of $\theta_a$;
\item for $a,b\in W$ define $\theta_{ab^{-1}}:X_{ba^{-1}}\rightarrow X_{ab^{-1}}$ by 

$$\theta_{ab^{-1}}(x)=\left\{\begin{array}{ll}
ay & \text{ if } |x|=\infty \text{ and } x=by,\\
(a\alpha,v) & \text{ if } x=(b\alpha,v),\\
(a,v) & \text{ if } x=(b,v);
\end{array} \right .$$

\item for the neutral element $0\in \F$ define $\theta_0:X_0\rightarrow X_0$ as the identity map;

\item for all the other elements $c$ of $\F$ define $\theta_c:X_{c^{-1}}\rightarrow X_c$ as the empty map.

\end{itemize}

The above maps together with the subsets $X_t$ form a partial action of $\F$ on $X$, that is $(\{\theta_t\}_{t\in \F}, \{X_t\}_{t\in \F})$ is such that $X_0=X$, $\theta_0=Id_x$, $\theta_c(X_{c^{-1}}\cap X_t)=X_{ct}\cap X_c$ and $\theta_c\circ\theta_t=\theta_{ct}$ in $X_{t^{-1}}\cap X_{t^{-1}c^{-1}}$. This partial action induces a partial action on the level of the $R$-algebra of functions (with point-wise sum and product) $F(X)$. 
More precisely, let $D$ be the subalgebra of $F(X)$ generated by all the finite sums of all the finite products of the characteristic maps $\{1_{X_A}\}_{A\in \mathcal{G}^0}$, $\{1_{bA}\}_{b\in W, A\in \mathcal{G}^0}$ and $\{1_{X_c}\}_{c\in \F}$. Also define, for each $t\in \F$, the ideals $D_t$ of $D$, as being all the finite sums of finite products of the characteristic maps $\{1_{X_t}1_{X_A}\}_{A\in \mathcal{G}^0}$, $\{1_{X_t}1_{bA}\}_{b\in W, A\in \mathcal{G}^0}$ and $\{1_{X_t}1_{X_c}\}_{c\in \F}$. Now, for each $c\in \F$, define the $R$-isomorphism $\beta_c:D_{c^{-1}}\rightarrow D_c$ by $\beta_c(f)=f\circ \theta_{c^{-1}}$. Then $(\{\beta_t\}_{t\in \F}, \{D_t\}_{t\in \F})$ is a partial action of $\F$ on $D$. 

\begin{remark}\label{spanDt} From now on we will use the notation $1_A$, $1_{bA}$ and $1_t$ instead of $1_{X_A}$, $1_{X_{bA}}$ and $1_{X_t}$, for $A\in \mathcal{G}^0$, $b\in W$ and $t\in \F$. Also, we have the following description of the ideals $D_t$:
$$D=D_0=\text{span}\{1_A, 1_c, 1_{bA}:A\in \mathcal{G}^0, c\in \F\setminus\{0\}, b\in W\}, $$ and, for each $t\in \F$, 
$$D_t=\text{span}\{1_t1_A, 1_t1_c, 1_t1_{bA}:A\in \mathcal{G}^0, c\in \F, b\in W\}. $$ 
\end{remark}

The key result in \cite{crossedproduct} that we need is the following theorem.

\begin{teorema}\cite{crossedproduct}\label{isom1} Let $\mathcal{G}$ be a countable ultragraph, $R$ be a unital commutative ring, and  $L_R(\mathcal{G})$ be the Leavitt path algebra of $\mathcal{G}$. Then there exists an $R$-isomorphism $\phi:L_R(\mathcal{G})\rightarrow D\rtimes_\beta \F$ such that $\phi(p_A)=1_A\delta_0$, $\phi(s_e^*)=1_{{e^{-1}}}\delta_{e^{-1}}$ and $\phi(s_e)=1_e\delta_e$, for each $A\in \mathcal{G}^0$ and $e\in \mathcal{G}^1$.
\end{teorema}

\begin{remark}
To prove Theorem~\ref{isom1} the authors first defined, using the universality of $L_R(\GG)$, a surjective homomorphism $\phi:L_R(\GG)\rightarrow D\rtimes_\beta\F$ such that $\phi(p_A)=1_A\delta_0$, $\phi(s_e)=1_e\delta_e$, and $\phi(s_e^{*})=1_{e^{-1}}\delta_{e^{-1}}$, for each $A\in \mathcal{G}^0$ and $e\in \mathcal{G}^1$. Up to this part of the proof the assumption on the cardinalilty of the ultragraph was not used. 
In particular, since $1_A\delta_0, 1_e\delta_e, 1_{e^{-1}}\delta_{e^{-1}}$ are all nonzero then $p_A, s_e$ and $s_{e^*}$ are all nonzero in $L_R(\GG)$ (for $\GG$ an arbitrary ultragraph). The injectivity of $\phi$ then followed by Theorem 3.2 of \cite{leavittultragraph}, which in turn relied on the Graded Uniqueness Theorem for Leavitt path algebras of countable graphs (Theorem 5.3 of \cite{Tomcommutativering}). However,  Theorem 5.3 of \cite{Tomcommutativering} also holds for arbitrary graphs, with the same proof. So, the injectivity of $\phi:L_R(G)\rightarrow D\rtimes_\beta\F$ also holds for arbitrary $\GG$, and we get the following theorem:
\end{remark}

\begin{teorema}\label{isom} Let $\mathcal{G}$ be an arbitrary ultragraph, $R$ be a unital commutative ring, and $L_R(\mathcal{G})$ be the Leavitt path algebra of $\mathcal{G}$. Then there exists an $R$-isomorphism $\phi:L_R(\mathcal{G})\rightarrow D\rtimes_\beta \F$ such that $\phi(p_A)=1_A\delta_0$, $\phi(s_e^*)=1_{{e^{-1}}}\delta_{e^{-1}}$ and $\phi(s_e)=1_e\delta_e$ for each $A\in \mathcal{G}^0$ and $e\in \mathcal{G}^1$.
\end{teorema}

\section{The reduction theorem}

The reduction theorem for Leavitt path algebras, see \cite{AAS, Malaga}, is an extremely useful tool in estabilishing various ring-theoretic properties of Leavitt path algebras (for example, the uniqueness theorems for Leavitt path algebras follow with mild effort from the reduction theorem). A version for relative Cohn path algebras was given in \cite{CG}, where it was also used as an important tool in the study of representations of these algebras. In our context the reduction theorem allow us to characterize faithful representations of Leavitt ultragraph path algebras arising from branching systems, but we expect it will also have applications in further studies of ultragraph Leavitt path algebras (for example, in Corollary~\ref{semiprime} we show that ultragraph Leavitt path algebras are semiprime). We present the theorem below, but first we recall the following:

\begin{definicao}[\cite{Tomsimplicity}] Let $\GG$ be an ultragraph. A closed path is a path
$\alpha \in \GG^*$ with $| \alpha | \geq 1$ and $s(\alpha)
\in r(\alpha)$. A closed path $\alpha$ is a cycle if $s(\alpha_i)\neq s(\alpha_j)$ for each $i\neq j$. An exit for a closed path is either an edge $e \in
\GG^1$ such that there exists an $i$ for which $s(e) \in
r(\alpha_{i})$ but $e \neq \alpha_{i+1}$, or a sink $w$ such that $w \in r(\alpha_i)$ for some $i$. We say that the ultragraph $\GG$ satisfies Condition~(L) if every closed path in $\GG$ has an
exit.
\end{definicao}

\begin{teorema}\label{reduction} Let $\GG$ be an arbitrary ultragraph, $R$ be a unital commutative ring and $0\neq x\in L_R(\GG)$. Then there are elements $\mu=\mu_1...\mu_n$ and $\nu=\nu_1...\nu_m\in L_R(\GG)$, with $\mu_i, \nu_j\in \GG^1\cup (\GG^1)^*$ for each $i$ and $j$, such that $0\neq \mu x\nu$ and either $\mu x\nu=\lambda p_A$, for some $A\in \GG^0$, or $\mu x\nu=\sum\limits_{i=1}^k\lambda_i s_{c}^i$, where $c$ is a cycle without exit.
\end{teorema}

\demo 

Recall that, by Theorem~\ref{isom}, there is an isomorphism $\Phi:L_R(\GG): D\rtimes_\beta \F$  such that $\Phi(p_A)=1_A\delta_0$ and $\Phi(s_ap_As_b^*)=1_{ab^{-1}}1_{aA}\delta_{ab^{-1}}$, where $a,b$ are paths in $\GG$ and $A\in \GG^0$. Moreover, $D\rtimes_\beta \F$ has a natural grading over $\F$, that is, $D\rtimes_\beta \F=\bigoplus_{t\in \F}D_t\delta_t$.

Let $x$ be a non-zero element in $L_R(\GG)$. We divide the proof in a few steps. The first one is the following. 

{\it Claim 1: There is a vertex $v\in G^0$ such that $xp_v\neq 0$}

First we prove that for each $t\in \F$, and $0\neq f_t\delta_t\in D_t\delta_t$, there is a vertex $v$ such that $f_t\delta_t1_v\delta_0\neq 0$.

Indeed, notice that for a non-zero element $f\delta_0\in D_0\delta_0$, since $f\neq 0$, there exists an element $y\in X$ such that $f(y)\neq 0$. Let $v=s(y)$. Then $(f1_v)(y)\neq 0$ and therefore $f\delta_01_v\delta_0=f1_v\delta_0\neq 0$. Similarly, for a given element $0\neq f\delta_{ab^{-1}}\in D_{ab^{-1}}\delta_{ab^{-1}}$, where $a,b$ are paths in $\GG$ (possibly one of them with length zero) notice that $\beta_{ba^{-1}}(f)\neq 0$ and hence there is an element $y\in X_{ba^{-1}}$ such that $(\beta_{ba^{-1}}(f))(y)\neq 0$. Let $v=s(y)$. Then $(\beta_{ba^{-1}}(f)1_v)(y)=(\beta_{ba^{-1}}(f))(y)\neq 0$, and it follows that $f\delta_{ab^{-1}}1_v\delta_0=\beta_{ab^{-1}}(\beta_{ba^{-1}}(f)1_v)\delta_{ab^{-1}}\neq 0$.

Now, since  $x \neq 0$, we have that $0\neq\Phi(x)=\sum\limits_{i=1}^n f_{t_i}\delta_{t_i}$, with $t_i\neq t_j$, and $f_{t_i}\delta_{t_i}\neq 0$ for each $i$. Fix some $i_0$ and chose a vertex $v$ such that $f_{t_{i_0}}\delta_{t_{i_0}}1_v\delta_0\neq 0$. Since $D\rtimes_\beta \F$ is $\F$-graded then $z=\sum\limits_{i=1}^nf_{t_i}\delta_{t_i}1_v\delta_0\neq 0$. Hence, applying $\Phi^{-1}$ to $z$, we get that $xp_v\neq 0$. 

{\it Claim 2: For each non-zero $x\in L_R(\GG)$ there exists an $y\in L_R(\GG)$ of the form $y=y_1...y_n$, with $y_i\in \GG^1$, such that $x y\neq 0$ and $xy$ has no ghost edges in its composition (a ghost edge is an element of $(\GG^1)^*$)}, that is, we can write $xy = \sum \lambda_j s_{a_j}p_{A_j}$ where $a_j$ are paths in the ultragraph $\GG$ and $A_j\in \GG^0$.

Write $x=\sum \lambda_i s_{a_i}p_{A_i}s_{b_i}^*$ with each $\lambda_i s_{a_i}p_{A_i}s_{b_i}^*\neq 0$ , where $a_i,b_i$ are paths in the ultragraph $\GG$, $A_i\in \GG^0$, and $\lambda_i\in R$ (see \cite[Theorem 2.5]{leavittultragraph}). If $b_i$ has length zero for each $i$ then we are done. Suppose that there exists an index $i$ such that $b_i$ has positive length. Write $b_i=eb$, where $e$ is an edge. Since $\lambda_i s_{a_i}p_{A_i}s_{b_i}^*\neq 0$ we have that $$0\neq \Phi(\lambda_i s_{a_i}p_{A_i}s_{b_i}^*)=\lambda_i 1_{a_i}1_{a_iA_i}1_{a_ib_i^{-1}}\delta_{a_ib_i^{-1}},$$ and so $\lambda_i 1_{a_i}1_{a_iA_i}1_{a_ib_i^{-1}}\neq 0$. Therefore we obtain that $$\Phi(\lambda_i s_{a_i}p_{A_i}s_{b_i}^*)\Phi(s_e)=\lambda_i1_{a_i}1_{a_iA_i}1_{a_ib_i^{-1}}\delta_{a_ib^{-1}}\neq 0.$$ Since $D\rtimes_\beta \F$ is $\F$-graded then $\Phi(x)\Phi(s_e)\neq 0$, and so $xs_e\neq 0$. Notice that $xs_e$ has less ghost edges in its composition than $x$. Repeating these arguments a finite number of times we obtain the conclusion of Claim 2.

{\it Claim 3. For each $0\neq x\in L_R(\GG)$ which is a sum of elements without ghost edges in its composition there are elements $y,z\in L_R(\GG)$, where $y=y_1...y_n$, $z=z_1...z_q$, and $y_i,z_j\in \GG^1\cup (\GG^1)^*\cup G^0$,  such that $yxz\neq 0$ and either $yxz=\lambda p_A$ for some $\lambda \in R$ and $A\in \GG^0$, or $yxz=\sum\limits_{i=1}^n\lambda_i s_c^i$ where $\lambda_i\in R$ and $c$ is a cycle without exit.} 

Let $0\neq x\in L_R(\GG)$ be any element without ghost edges in its composition. Write $x=\sum\limits_{j\in M}\beta_jp_{A_j}+\sum\limits_{i\in N} \lambda_is_{a_i}$, where $A_j\in \GG^0$, $a_i$ are paths with  positive length, $\beta_j,\lambda_i$ are nonzero elements in $R$, and the number of summands describing $x$ is the least possible. Notice that $M$ or $N$ could be empty. Define $m$ as the cardinality of $M$ and $n$ as the cardinality of $N$.

We prove the claim using an induction argument over the (minimal) number of summands in $x\in L_R(\GG)$  without ghost edges in its composition.

If $m+n=1$ then $x=\beta_1P_{A_1}$, or $x=\lambda_1s_{a_1}$ (in which case  $s_{a_1}^*x=\lambda_1p_{r(a_1)}$), and so we are done.

Now, suppose the induction hypothesis holds and let $x=\sum\limits_{j\in M}\beta_jp_{A_j}+\sum\limits_{i\in N} \lambda_is_{a_i}$ (with minimal number of summands). We prove below that Claim~3 holds for this $x$.

Suppose that $N$ is empty. By Claim~1 there is a vertex $v$ such that $0\neq xp_v$. Hence $xp_v=(\sum\limits_{j\in M:v\in A_j}\beta_j)p_v$, and we are done.

Now, suppose that $N$ is nonempty, say $N=\{1,2,...,n\}$. Moreover assume, without loss of generality, that $|a_i|\leq |a_{i+1}|$ for each $i\in \{1,...,n-1\}$. 

By Claim~1 there exist a vertex $v$ such that $xp_v\neq 0$. If $m>1$ then $$xp_v=\sum\limits_{j\in M}\beta_j1_{A_j}p_v + \sum\limits_{i\in N}\lambda_is_{a_i}p_v = (\sum\limits_{j:v\in A_j}\beta_j)p_v+\sum\limits_{i\in N}\lambda_is_{a_i}p_v,$$
 which has less summands than $x$, and so we may apply the induction hypothesis on $xp_v$. The same holds if $s_{a_i}p_v=0$ for some $i$.
Therefore we are left with the cases when $m=0$ or $m=1$, and $s_{a_i}p_v\neq 0$ for each $a_i$.

Before we proceed, notice that each $a_i$ in the description of $x$ is of the form $a_i=a_i^1...a_i^{k_i}$, where $a_i^j$ are edges.
If there are $i,k$ such that $a_i^j\neq a_k^j$ then $s_{a_i}^*s_{a_k}=0$, and $s_{a_i}^*x\neq 0$ (since $L_K(\GG)$ is $\F$ graded and $s_{a_i}^* s_{a_i}=p_{r(a_i)}\neq 0$). Therefore $s_{a_i}^*x$ has less summands then $x$, possibly including ghost edges. Applying Claim~2 to $s_{a_i}^*x$, and then applying the induction hypothesis, we obtain the desired result for $x$. So, we may suppose that each path $a_i$ is the beginning of the path $a_{i+1}$.

Recall that to finish the proof we need to deal with two cases:  $$xp_v=\beta p_v+\sum\limits_{i=1}^n\lambda_i s_{a_i}p_v \text{, or }xp_v=\sum\limits_{i=1}^n\lambda_i s_{a_i}p_v,$$
 where $|a_i|\leq |a_{i+1}|$ and $a_i$ is the beginning of $a_{i+1}$ for each $i$.
 
If  $xp_v=\sum\limits_{\lambda_i}s_{a_i}p_v$, then $s_{a_1}^*xp_v=\lambda_1p_v+\sum\limits_{i=2}^ns_{a_1}^*s_{a_i}p_v$. Therefore it is enough to deal with the case $xp_v=\beta p_v+\sum\limits_{i=1}^n\lambda_i s_{a_i}p_v$. 

Notice that (by the $\F$-grading on $L_R(\GG)$) $p_vxp_v\neq 0$, and $p_vxp_v=\beta p_v+\sum\limits_{i=1}^n \lambda_ip_vs_{a_i}p_v$. If there is a $j$ such that $s(a_j)\neq v$ then $p_vs_{a_j}=0$, and so we may apply the induction hypothesis on $p_vxp_v$. Therefore we are left with the case when $v=s(a_i)$ for each $i$, what implies that each $a_i$ is a closed path based on $v$.

Let $c=a_1$, and write $c=c_1...c_k$. If $c$ has an exit, then either there exists an edge $e\neq c_{j+1}$ such that $s(e)\in r(c_j)$, or there exists a sink $w$ in $r(c_i)$ for some $i$. In the first case, notice that, since each $a_i$ is of the form $a_i=c\overline{a_i}$, then $s_{c_1...c_je}^*s_{a_i}=0$ for each $i$. Then $$s_{c_1...c_je}^*xp_vs_{c_1...c_je}=s_{c_1...c_je}^*s_{c_1...c_je}=\beta p_{r(e)}\neq 0.$$
Now suppose that there is a sink $w$ in $r(c_q)$ for some $q$. Then for each $i$, $p_ws_{c_1...c_q}^*s_{a_i}s_{c_1...c_q}=p_ws_{a_i^{k+1}}...s_{a_i^{k_i}}s_{c_1...c_q}=0$ since $w$ is not the source of any edge. Then $p_ws_{c_1...c_q}^*xp_vs_{c_1...c_q}=\beta p_w p_{r(c_q)}=\beta p_w\neq 0$.

So it remains the case when $c$ has no exit. Since $c$ is based on $v$ and has no exit, and each $a_i$, for $i>1$, is a closed path based on $v$ then $a_i$ equals $c^j$ for some $j$. Notice that since $c$ is a closed path without exit then $c=d^n$, where $d$ is a cycle without exit and $n\in \N$. Therefore $s_{a_1}xp_v$ has the desired form.
\fim

As a first example of potential applications of the reductions theorem we show below that ultragraph Leavitt path algebras are semiprime (for Leavitt path algebras of graphs over fields this is Proposition~2.3.1 in \cite{AAS}). Recall that a ring $R$ is said to be semiprime if, for every ideal $I$ of $R$, $I^2=0$ implies $I=0$.

\begin{corolario}\label{semiprime} Let $\GG$ be an arbitrary ultragraph and $R$ be a unital commutative ring with no zero divisors. Then the ultragraph Leavitt path algebra $L_R(\GG)$ is semiprime.
\end{corolario}

\demo
Let $I$ be a nonzero ideal and $x$ in $I$ be nonzero. By the Reduction Theorem~\ref{reduction} there are elements $\mu=\mu_1...\mu_n$ and $\nu=\nu_1...\nu_m\in L_R(\GG)$, with $\mu_i, \nu_j\in \GG^1\cup (\GG^1)^*$ for each $i$ and $j$, such that $0\neq \mu x\nu$ and either $\mu x\nu=\lambda p_A$, for some $A\in \GG^0$, or $\mu x\nu=\sum\limits_{i=1}^k\lambda_i s_{c}^i$, where $c$ is a cycle without exit. Since $p_A$ is an idempotent if $\mu x\nu=\lambda p_A$ then we are done. If $\mu x\nu=\sum\limits_{i=1}^k\lambda_i s_{c}^i$ then, by the $\F$ grading (see Theorem~\ref{isom}), we have that $(\mu x\nu)^2\neq 0$ and we are done.
\fim

\section{Algebraic branching systems and the induced representations}

In this section we start the study of representations of ultragraph Leavitt path algebras via branching systems. Motivated by the relations that define an ultragraph Leavitt path algebra we get the following definition. 

\begin{definicao}\label{branchsystem}
Let $\mathcal{G}$ be an ultragraph, $X$ be a set and let $\{R_e,D_A\}_{e\in \mathcal{G}^1,A\in \mathcal{G}^0}$ be a family of subsets of $X$. Suppose that
\begin{enumerate}
\item\label{R_e cap R_f =emptyset if e neq f} $R_e\cap R_f =\emptyset$, if $e \neq f \in \mathcal{G}^1$;
\item $D_\emptyset=\emptyset, \ D_A \cap D_B= D_{A \cap B}, \text{ and } D_A \cup D_B= D_{A \cup B}$ for all $A, B \in \mathcal{G}^0$;
\item $R_e\subseteq D_{s(e)}$ for all $e\in \mathcal{G}^1$;
\item\label{D_v=cup_{e in s^{-1}(v)}R_e} $D_v= \bigcup\limits_{e \in s^{-1}(v)}R_e$, if $0 <\vert s^{-1}(v) \vert<\infty$; and
\item\label{funcoes} for each $e\in \mathcal{G}^1$, there exist two bijective maps $f_e:D_{r(e)}\rightarrow R_e$ and $f_e^{-1}:R_e \rightarrow D_{r(e)}$ such that $f_e\circ f_e^{-1}=\id_{R_e}$ and, $f_e^{-1}\circ f_e=\id_{D_{r(e)}}$.
\end{enumerate}

We call $\{R_e,D_A,f_e\}_{e \in \mathcal{G}^1,A \in \mathcal{G}^0}$ a $\mathcal{G}$-algebraic branching system on $X$ or, shortly, a $\GG$-branching system.
\end{definicao}

Given any ultragraph $\GG$ we deal with the existence of branching systems associated to $\GG$ in the proposition below. 

\begin{proposicao}\label{existenceofabranchingsystem}
Let $\mathcal{G}$ be an ultragraph such that $s^{-1}(v)$ and $r(e)$ are finite or countable for each vertex $v\in G^0$ and each edge $e$. Then there exists a $\mathcal{G}$-branching system.
\end{proposicao}

\demo Let $X=[0,1)\times (E^1\cup G^0)$. For each edge $e$ define $R_e=[0,1)\times \{e\}$. For each sink $u\in G^0$ define $D_u=[0,1)\times \{u\}$, and for each non sink $v\in G^0$ let $D_v=\bigcup\limits_{e\in s^{-1}(v)}R_e=[0,1)\times s^{-1}(v)$.  Define $D_\emptyset=\emptyset$ and, for each non-empty $A \in \mathcal{G}^0$, let $$D_A:=\bigcup_{v \in A}D_v.$$

It is easy to see that $\{R_e,D_A\}_{e \in \mathcal{G}^1,A \in \mathcal{G}^0}$ satisfies Condition~(\ref{R_e cap R_f =emptyset if e neq f})--(\ref{D_v=cup_{e in s^{-1}(v)}R_e}) of Definition~\ref{branchsystem}. We prove that Condition~\ref{funcoes} is satisfied.

Fix $e \in \mathcal{G}^1$. We have to construct $f_e$ and $f_e^{-1}$ that satisfy Condition~\ref{funcoes}. Since $D_{r(e)}=\bigcup\limits_{v\in r(e)}D_v$, and $r(e)$ and $s^{-1}(v)$ are countable sets, then $D_{r(e)}=[0,1)\times J$ where $J$ is a finite or countable set.
If $J$ is finite then let $J=\{c_1,...,c_n\}$, and for each $i\in \{1,...,n\}$ define $F_i:[\frac{i-1}{n},\frac{i}{n})\times \{e\}\to [0,1)\times \{c_i\}$ by $F_i(x,e)=(nx-i+1,c_i)$, which is bijective. Piecing together $F_i$'s yields $f_e$, and piecing together $F_i^{-1}$'s yields $f_e^{-1}$. If $J$ is infinite let $J=\{c_i\}_{i\in \N}$, and for each $i\in \N$ define $F_i:[1-\frac{1}{i}, 1-\frac{1}{i+1})\times\{e\}\to [0,1)\times \{c_i\}$ by $F_i(x,e)=((i+1)ix-(i+1)(i-1),c_i)$, which is bijective. Again, 
piecing together $F_i$'s yields $f_e^{-1}$, and piecing together $F_i^{-1}$'s yields $f_e$.
\fim
\begin{remark} The result above extends Theorem~3.1 in \cite{represleavitt} to uncountable graphs.
\end{remark}

Next we describe how to construct representations of ultragraph Leavitt path algebras from branching systems. 

Let $\mathcal{G}$ be an ultragraph, $R$ be a commutative unital ring and $X$ be a $\GG$-algebraic branching system. Denote by $M$ the set of all maps from $M$ to $R$. Notice that $M$ is an $R$-module, with the usual operations. For a subset $Y\subseteq X$ we write $1_Y$ to denote the   characteristic function of $Y$, that is $1_Y(x)=1$ if $x\in Y$, and $1_Y(x)=0$ otherwise (clearly $1_Y$ is an element of $M$). Finally, we denote by $Hom(M)$ the set of all $R$-homomorphisms from the $R$-module $M$ to $M$.

For each $e\in \mathcal{G}^1$ and $\phi\in M$ define 

$$S_e(\phi)=\phi\circ f_{e}^{-1}\cdot1_{R_e} \text{ and }S_e^*(\phi)=\phi\circ f_{e}\cdot1_{D_{r(e)}},$$ and for each $A\in \mathcal{G}^0$ and $\phi\in M$ define $$P_A(\phi)=1_{D_A}\phi.$$

Clearly, for each $e\in \GG^1$ and $A\in \GG^0$, we have that $S_e, S_e^*, \text{ and } P_A\in Hom(M)$.

\begin{remark} To simplify our notation, for each $e\in \mathcal{G}^1$ and $\phi\in M$, we write $S_e(\phi)=\phi\circ f_{e}^{-1}$ and $S_e^*(\phi)=\phi\circ f_e$ instead of $S_e(\phi)=\phi\circ f_{e}^{-1}\cdot 1_{R_e}$ and $S_e^*(\phi)=\phi\circ f_e \cdot 1_{D_{r(e)}}$.
\end{remark}

We end the section describing the representations induced by branching systems associated to an ultragraph $\GG$.

\begin{proposicao}\label{repinducedbybranchingsystems}
Let $\mathcal{G}$ be an ultragraph and let $\{R_e,D_A,f_e\}_{e\in \mathcal{G}^1, A\in \mathcal{G}^0}$ be a $\GG$-algebraic branching system on $X$. Then there exists a unique representation $\pi:L_R(\GG) \rightarrow Hom(M)$ such that $\pi(s_e)(\phi)= S_e$, $\pi(s_e^*)=S_e^*$, and $\phi(p_A)= P_A$, or equivalently, such that $\pi(s_e)(\phi)=\phi\circ f_{e}^{-1}$, $\pi(s_e^*)(\phi)=\phi\circ f_e$ and $\phi(p_A)(\phi)=\phi 1_{D_A}$, for each $e\in \GG^1$, $A\in \GG^0$, and $\phi\in M$.
\end{proposicao}

\demo 

By the universal property of ultragraph Leavitt path algebras we only need to verify that the family $\{S_e, S_e^*, P_A\}_{e\in \GG^1, A \in \GG^0}$ satisfy relations 1 to 4 in Definition~\ref{def of ultragraph algebra}. We show below how to verify relation 4 and leave the others to the reader. 

Note that for $e\in \GG^1$ and $\phi\in M$ we get $\pi(s_e)\pi(s_e^*)(\phi)=\pi(s_e)(\phi\circ f_e)=\phi\circ f_e\circ f_e^{-1}=\phi \cdot 1_{R_e}$. Now, let $v$ be a vertex such that $0<|s^{-1}(v)|<\infty$. Then for each $\phi\in M$,  we have that $$\sum\limits_{e\in s^{-1}(v)}\pi(s_e)\pi(s_e^*)(\phi)=\sum\limits_{e\in s^{-1}(v)}\phi \cdot 1_{R_e}=\phi (1_{\bigcup\limits_{e\in s^{-1}(v)}R_e})=\phi\cdot 1_{D_v}=\pi(p_v)(\phi),$$ and relation 4 in Definition~\ref{def of ultragraph algebra} is proved.
\fim

\begin{remark}\label{finitesupport} In the previous theorem we can also take $M$ as being the $R$-module of all the maps from $X$ to $R$ with finite support, instead of all the maps from $X$ to $R$. 
\end{remark}

\section{Faithful Representations of ultragraph Leavitt path algebras via Branching Systems}

In this section we describe faithfulness of the representations induced by branching systems in terms of dynamical properties of the branching systems, and in terms of combinatorial properties of the ultragraph. Our first result follows below, linking faithfulness of the representations with the dynamics of the branching systems.

\begin{teorema}\label{faithfulrep1}
Let $\GG$ be an ultragraph, and $\{R_e,D_A,f_e\}_{e\in \mathcal{G}^1, A\in \mathcal{G}^0}$ be an branching system such that $D_A\neq \emptyset$ for each $\emptyset \neq A\in \GG^0$. Then the induced representation of $L_R(\GG)$ (from Proposition~\ref{repinducedbybranchingsystems}) is faithful if, and only if, for each closed cycle without exit $c=c_1...c_k$, and for each set finite set $F\subseteq \N$,  there exists an element $z_0\in D_{r(c)}$ such that $ f_c^n(z_0)\neq z_0$, for each $n\in F$ (where $f_c=f_{c_1}...f_{c_n}$).
\end{teorema}

\demo Let $\pi:L_R(\GG)\rightarrow Hom(M)$ be the homomorphism induced by Proposition~\ref{repinducedbybranchingsystems}. 

Let $0\neq x\in L_R(\GG)$. By Theorem \ref{reduction} there are elements $y,z\in L_R(\GG)$ such that 
$0\neq y x z$ and either $y x z=\lambda p_A$, for some $A\in \GG^0$, or $y x z=\sum\limits_{i=1}^k\lambda_i s_{c}^i$, where $c$ is a cycle without exit.

If $yxz=\lambda p_A$, for some $A\in \GG^0$, then $\pi(x)\neq 0$ since 
$$\pi(y)\pi(x)\pi(z)(1_{D_A})=\pi(yxz)(1_{D_A})=\pi(\lambda p_A)(1_{D_A})=\lambda 1_{D_A}\neq 0.$$

Suppose that $yxz=\sum\limits_{i=1}^m\lambda_i s_{c}^i$, where $c$ is a cycle without exit. Let $j$ be the least of the elements $i$ such that $\lambda_i\neq 0$. Define $\mu=(s_c^j)^*yxz$, and note that $\mu=\lambda_j p_v+\sum\limits_{i=1}^{m-j} \widetilde{\lambda_i}s_{c^i}$, where $\widetilde{\lambda_i}=\lambda_{j+i}$ and $v=r(c)$.

Let $z_0\in D_v=D_{r(c)}$ be such that $f_c^i(z_0)\neq z_0$ for each $i\in \{1,...,m-j\}$ (such $z_0$ exists by hypothesis). Let $\delta_{z_0}\in M$ be the map defined by $\delta_{z_0}(x)=1$ if $x=z_0$, and $\delta_{z_0}(x)=0$ otherwise. Notice that, for each $i\in \{1,...,m-j\}$, we have
$$\pi(s_c^i)(\delta_{z_0})(z_0)=\delta_{f_c^i(z_0)}(z_0)=0.$$ Therefore $$\pi(\mu)(\delta_{z_0})(z_0)=\lambda_j\pi(p_v)(\delta_{z_0})(z_0)=\lambda_j\neq 0,$$ from where $\pi(\mu)\neq 0$ and hence $\pi(x)\neq 0$ (since $\mu=(s_c^j)^*yxz$). 

For the converse suppose that there exist a $j_0$, and a cycle $c$ without exits based at $w$ such that $f_c^{j_0}(z)=z$ for every $z \in D_w$. Then we have that $\pi(s_{c^{j_0}})=\pi(p_w)$. To see that $p_w\neq s_{c^{j_0}}$ use Theorem \ref{isom} and the $\F$-grading of $D\rtimes_\beta\F$. So, $\pi$ is not injective.


\fim

For ultragraphs that satisfy Condition~(L) the above theorem has a simplified version. Recall that an ultragraph $\GG$ satisfies Condition~(L) if every closed path in $\GG$ has an
exit.

\begin{corolario}\label{faithfulrep}
Let $\GG$ be an ultragraph satisfying Condition~$(L)$, and let $\{R_e,D_A,f_e\}_{e\in \mathcal{G}^1, A\in \mathcal{G}^0}$ be a branching system with $D_A\neq \emptyset$, for each $\emptyset \neq A\in \GG^0$. Then the induced representation of $L_R(\GG)$ (from Theorem~ \ref{repinducedbybranchingsystems}) is faithful.
\end{corolario}

For ultragraphs such that $s^{-1}(v)$ and $r(e)$ are finite or countable (for each vertex $v$ and each edge $e$) we also get the converse of the previous corollary. This is our next result.

\begin{teorema}
Let $\GG$ be an ultragraph such that $s^{-1}(v)$ and $r(e)$ are finite or countable for each $v\in G^0$ and $e\in E^1$. Then $\GG$ satisfies Condition~$(L)$ if, and only if, for every algebraic branching system $\{R_e,D_A,f_e\}_{e\in \mathcal{G}^1, A\in \mathcal{G}^0}$ with $D_A\neq \emptyset$, for each $\emptyset \neq A\in \GG^0$, the induced representation of $L_R(\GG)$ (from Proposition~\ref{repinducedbybranchingsystems}) is faithful.
\end{teorema}

\demo By Corollary~\ref{faithfulrep} we only need to prove the converse.
We do this proving the contrapositive, i.e., we prove that if $\GG$ does not satisfy condition $(L)$ then there exists an algebraic branching system $\{R_e,D_A,f_e\}_{e\in \mathcal{G}^1, A\in \mathcal{G}^0}$ with $D_A\neq \emptyset$, for each non-empty$A\in \GG^0$, such that the representation induced by Proposition~\ref{repinducedbybranchingsystems} is not faithful.

Suppose that $\GG$ does not satisfy condition $(L)$. Then there exists a cycle $\alpha=(\alpha_1,\dots,\alpha_n)$ such that $\vert r(\alpha_i) \vert=1$ and $s^{-1}(s(\alpha_i))=\{\alpha_i\}$, for all $i=1,\dots,n$, and $\alpha_i \neq \alpha_j$ if $i \neq j$. 

By the proof of Theorem~\ref{existenceofabranchingsystem} there is a $\mathcal{G}$-branching system on $[0,1)\times (G^0\cup E^1)$ such that $D_{r(\alpha_{i-1})}=D_{s(\alpha_i)}=R_{\alpha_i}=[0,1)\times \{\alpha_i\}$, and $f_{\alpha_i}:D_{r(\alpha_i)}\to R_{\alpha_i}$ is the bijective affine map $f_{\alpha_i}(x,\alpha_{i+1})=(x,\alpha_i)$, for each $i=1,\dots,n$. Notice that $f_{\alpha_1}\circ f_{\alpha_2}...\circ f_{\alpha_n}=\id_{R_{\alpha_1}}$. Let $\pi$ be the representation associated with this branching system (as in Proposition~\ref{repinducedbybranchingsystems}). Then, for each $\phi\in M$, we have $$\pi(s_{\alpha_n}^*)...\pi(s_{\alpha_1}^*)(\phi)=\phi\circ(f_{\alpha_1}\circ...\circ f_{\alpha_n}).1_{D_{s(\alpha_1)}}=\phi.1_{D_s{(\alpha_1)}}=\pi(p_{s(\alpha_1)})(\phi),$$ and hence $\pi(s_{\alpha_n}^*)...\pi(s_{\alpha_1}^*)=\pi(p_{s(\alpha_1)})$.

To see that $\pi$ is not faithful it remains to show that $s_{\alpha_n} ^*...s_{\alpha_1}^*\neq p_{s(\alpha_1)}$. For this, we construct a branching system as follows: Let $\{D_A\}_{A\in \GG^0}$ and $\{R_e\}_{e\in \GG^1}$ be as above, and chose maps $\widetilde{f}_{\alpha_1},...\widetilde{f}_{\alpha_n}$ such that $\widetilde{f}_{\alpha_1}\circ...\circ \widetilde{f}_{\alpha_n}\neq Id_{R_{\alpha_1}}$. Let $x_0\in R_{\alpha_1}$ be such that $(f_{\alpha_1}\circ...\circ f_{\alpha_n})(x_0)\neq x_0$, and chose an element $\varphi\in M$ such that $\varphi(x_0)=1$ and $\varphi\circ(f_{\alpha_1}\circ ...\circ f_{\alpha_n})(x_0)=0$.
Let $\widetilde{\pi}$ be the representation of $L_R(\GG)$ obtained by Proposition~\ref{repinducedbybranchingsystems} from this branching system. Then $$\widetilde{\pi}(s_{\alpha_n}^*\ldots s_{\alpha_1}^*)(\varphi)(x_0)=\varphi(\widetilde{f}_{\alpha_1}\circ...\circ\widetilde{f}_{\alpha_n})(x_0)=0\neq 1=\varphi(x_0)=\widetilde{\pi}(p_{s_{\alpha_1}})(\varphi)(x_0),$$ and so $\widetilde{\pi}(s_{\alpha_n}^*)...\widetilde{\pi}(s_{\alpha_1}^*)\neq \widetilde{\pi}(p_{s(\alpha_1)})$. It follows that $s_{\alpha_n} ^*...s_{\alpha_1}^*\neq p_{s(\alpha_1)}$, and hence the representation $\pi$ is not faithful.
\fim

We end this section showing that it is always possible to construct faithful representations of $L_R(\GG)$ arising from branching systems.

\begin{proposicao} Let $\GG$ be an ultragraph such that $s^{-1}(v)$ and $r(e)$ are finite or countable for each $v\in G^0$ and $e\in E^1$. Then there exist a $\GG$-branching system $\{R_e, D_A, f_e\}_{e\in \GG^1, A\in \GG^0}$ such that the representation induced from Proposition~\ref{repinducedbybranchingsystems} is faithful.
\end{proposicao}

\demo 
Let $\{R_e,D_A,f_e\}_{e\in \GG^1, A \in \GG^0}$ be the $\GG$-branching system obtained in Theorem \ref{existenceofabranchingsystem}. Following Theorem \ref{faithfulrep1}, all we need to do is to redefine some maps $f_e:D_{r(e)}\rightarrow R_e$ to get a branching system such that for each closed cycle $c$ without exit, and for each finite set $F\subseteq \N$, there is an element $z_0\in D_{r(c)}$ such that $f_c^i(z_0)\neq z_0$ for each $i\in F$. Let $c=c_1...c_n$ be a closed cycle without exit, where $c_i$ are edges. Recall from the proof of Theorem~ \ref{existenceofabranchingsystem} that $R_e=[0,1)\times \{e\}$ for each edge $e$, and $D_v=\bigcup\limits_{e\in s^{-1}(v)}R_e=[0,1)\times s^{-1}(v)$ for each non sink $v$. Since $c=c_1...c_n$ is a closed cycle without exit then $D_{r(c_i)}=[0,1)\times c_{i+1}$ for each $i\in \{1,...,n-1\}$, and $D_{r(c_n)}=[0,1)\times c_1$. Now, let $\theta\in(0,1)$ be an irrational number. Define, for each $i\in \{1,...,n-1\}$, $f_{c_i}:[0,1)\times c_{i+1}\rightarrow [0,1)\times c_i$ by $f_{c_i}(x,c_{i+1})=((x+\theta)mod(1), c_i)$ and define $f_{c_n}:[0,1)\times c_1\rightarrow [0,1)\times c_n$ by $f_{c_n}(x,c_1)=((x+\theta)mod(1), c_n)$. Let $f_c:D_{r(c_n)}\rightarrow R_{c_1}$ be the composition $f_c=f_{c_1}...f_{c_n}$, and notice that for each rational number $x\in [0,1)$ we get $f_c(x,c_1)=(y,c_1)$, where $y$ is an irrational number. Therefore $f_c(x,c_1)\neq (x,c_1)$ for each rational number $x\in [0,1)$.
\fim

\section{Permutative representations}

Our aim in this section is to show that under a certain condition over an ultragraph $\GG$, each $R$-algebra homomorphism $\pi: L_R(\GG)\rightarrow A$ has a sub-representation associated to it which is equivalent to a representation induced by a $\GG$-algebraic branching system. 


Recall that given an $R$-algebra $A$ (where $R$ is a unital commutative ring), there exist an $R$-module $M$ and an injective $R$-algebra homomorphism $\varphi:A\rightarrow Hom_R(N)$ (see Section~5 of \cite{represleavitt} for example).  
Composing a homomorphism from $L_R(\GG)$ to an $R$-algebra $A$ with the previous homomorphism $\varphi$ we get a homomorphism from $L_R(\GG)$ to $Hom_R(N)$. So we will only consider representations from $L_R(\GG)$ with image in  $Hom_R(N)$, for some $R-$module $N$.

Next we set up notation and make a few remarks regarding a representation (an $R$-homomorphism) $\pi:L_R(\GG)\rightarrow Hom_R(N)$, where $N$ is an $R$-module.

For each $e\in \mathcal{G}^1$ define $M_e=\pi(s_es_e^*)(N)$, and for each $A\in \mathcal{G}^0$ define $N_A=\pi(p_A)(N)$. Notice that $N_e$ and $N_A$ are sub-modules of $N$ and:
\begin{itemize}
\item $N_e\cap N_f=\{0\}$ for each $e\neq f$;
\item $N_e\subseteq N_{s(e)}$ for each $e\in \mathcal{G}^1$;
\item for each non sink $v\in \mathcal{G}^0$ we have $N_v=\left(\bigoplus\limits_{s(e)=v}N_e\right)\bigoplus O_v$, where $O_v$ is a submodule of $N_v$. If $|s^{-1}(v)|<\infty$ then $O_v=\{0\}$;
\item for each $A\in \mathcal{G}^0$ $N_A=\left(\bigoplus\limits_{v\in A}N_v\right)\bigoplus O_A$, where $O_A$ is a submodule of $N_A$. If $A$ is finite then $O_A=\{0\}$. 
\item $\pi(s_e^*)_{|_{N_e}}:N_e\rightarrow N_{r(e)}$ is an isomorphism for each $e\in \mathcal{G}^1$, with inverse $\pi(s_e)_{|_{N_{r(e)}}}:N_{r(e)}\rightarrow N_e$.

\end{itemize}


Permutative representations were originally defined by Bratteli-Jorgensen in the context of representations of Cuntz algebras (see \cite{BJ}). Without knowledge of Bratteli-Jorgensen definition the authors defined permutative representations, originally in the context of graph C*-algebras and then for Leavitt path algebras and ultragraph C*-algebras (see \cite{unitaryequivalence, represleavitt, ultragraphrep}), as representations that satisfy Condition~(B2B). Next we define this condition in the context of Leavitt path algebra of ultragraphs.

\begin{definicao}\label{permutativerep} Let $\pi:L_R(\GG)\rightarrow Hom_R(N)$ be an $R$-homomorphism. We say that $\pi$ is permutative if there exist  bases $B$ of $N$, $B_v$ of $N_v$, $B_{r(e)}$ of $N_{r(e)}$ and $B_e$ of $N_e$, for each $v\in G^0$ and $e\in \GG^1$, such that:

\begin{enumerate}
\item $B_v\subseteq B$ and $B_{r(e)}\subseteq B$ for each $v\in G^0$ and $e\in \GG^1$,
\item $B_v\subseteq B_{r(e)}$ for each $v\in r(e)$,
\item $B_e\subseteq B_{s(e)}$ for each $e\in \GG^1$,
\item $\pi(s_e)(B_{r(e)})=B_e$ for each $e\in \GG^1$ $(B2B)$.
\end{enumerate}
\end{definicao}

\begin{remark} The last condition of the previous definition is equivalent to $\pi(s_e^*)(B_e)=B_{r(e)}$, since the map $\pi(s_e):N_{r(e)}\rightarrow N_e$ is an $R$-isomorphism, with inverse $\pi(s_e^*):N_e\rightarrow N_{r(e)}$. 
\end{remark}

Notice that associated to a permutative representation $\pi:L_R(\GG)\rightarrow Hom_R(N)$ there are basis of $N$, $N_v$, $N_{r(e)}$, and $N_e$, for every $v\in G^0$ and $e\in \GG^1$. From these basis we can build a basis of $N_A$, for every $A\in \GG^0$. More precisely, for each $A \in \GG^0$, following Lemma~\ref{concrete des of G^0}, write $A=\left(\bigcap\limits_{e\in X_1}r(e) \right)\cup...\cup\left(\bigcap\limits_{e\in X_n}r(e) \right)\bigcup F$,  and define $$B_A=\left(\bigcap\limits_{e\in X_1}B_{r(e)} \right)\cup...\cup\left(\bigcap\limits_{e\in X_n}B_{r(e)} \right)\bigcup\limits_{v\in F} B_v.$$ We then have the following.

\begin{lema}\label{welldefined} Let $\pi:L_R(\GG)\rightarrow Hom_R(N)$ be a permutative representation and $A \in \GG^0$. Then the set $B_A$ is well defined (that is, does not depends on the description of $A$) and is a basis of $N_A$. Moreover, if $D\in \GG^0$ then $B_{A\cup D}=B_A\cup B_D$ and $B_{A\cap D}=B_A\cap B_D$.
\end{lema}

\demo First we prove that $B_A$ is a basis of $N_A=\pi(p_A)(N)$. For each $h\in B_v$, with $v\in F$, we have $\pi(p_A)(\pi(p_v)(h))=\pi(p_v)(h)=h$ and so $h\in N_A$. Moreover, for $h\in \bigcap\limits_{e\in X_i}B_{r(e)}$, we get $\pi(p_A)(\pi(p_{\bigcap\limits_{e\in X_i}r(e)})(h))=\pi(p_{\bigcap\limits_{e\in X_i}r(e)})(h)=h$ and so $h\in N_A$. Therefore, $B_A\subseteq N_A$. Next we show that $span(B_A)=N_A$.
Notice that, by definition, $N_{r(e)}=span(B_{r(e)})$ for each edge $e$, and $N_v=span(B_v)$ for each vertex $v$. To show that $span(B_A)=N_A$ for any $A\in \GG^0$ we first need to prove the following:

{\it Claim 1: For each finite set $X\subseteq \GG^1$ it holds that $span(\bigcap\limits_{e\in X}B_{r(e)})=\bigcap\limits_{e\in X}span(B_{r(e)})$.}

The inclusion $span(\bigcap\limits_{e\in X}B_{r(e)})\subseteq\bigcap\limits_{e\in X}span(B_{r(e)})$ is obvious. To prove the other inclusion, let $e,f$ be two edges, and let $h\in span(B_{r(e)})\cap span(B_{r(f)})$. Write $h=\sum\limits_{i=1}^n\alpha_i h_i=\sum\limits_{j=1}^m\beta_j k_j$, with $h_i\in B_{r(e)}$ for each $i$, $k_j\in B_{r(f)}$ for each $j$, and all $\alpha_i$ and $\beta_j$ non-zero. Then $\sum\limits_{i=1}^n\alpha_ih_i-\sum\limits_{j=1}^m\beta_j k_j=0$. If $h_i \notin \{k_1,...,k_m\}$ for some $i$ then $\alpha_i=0$ (since $\{h_1,...,h_n\}\cup\{k_1,...,k_m\}\subseteq B$ and $B$ is linearly independent), which is impossible since $\alpha_i\neq 0$ for each $i$. So we get $\{h_1,...,h_n\}\subseteq \{k_1,...,k_m\}$. By the same arguments, applied to $k_i$, we obtain that $\{k_1,...,k_m\}\subseteq \{h_1,...,h_n\}$. Therefore $\{k_1,...,k_m\}= \{h_1,...,h_n\}$ and hence $\{h_1,...,h_n\}\subseteq B_{r(e)}\cap B_{r(f)}$ and $h\in span(B_{r(e)}\cap B_{r(f) })$. The claim now follows by inductive arguments over the the cardinality of $X$.

We now prove that $N_A=span(B_A)$

Suppose first that $A=\bigcap\limits_{e\in X_1}r(e)\cup\bigcap\limits_{f\in X_2}r(f)$. Then 
$$N_A=\pi(p_A)(N)=\prod\limits_{e\in X_1}\pi(p_{r(e)})(N)+\prod\limits_{f\in X_2}\pi(p_{r(f)})(N)-\prod\limits_{e\in X_1}\pi(p_{r(e)})\prod\limits_{f\in X_2}\pi(p_{r(f)})(N).$$ Notice that $\prod\limits_{e\in X_1}\pi(p_{r(e)})(N)\subseteq \pi(p_{r(e)})(N)=span(B_{r(e)})$ for each $e\in X_1$ and so $\prod\limits_{e\in X_1}\pi(p_{r(e)})(N)\subseteq \bigcap\limits_{e\in X_1}span(B_{r(e)})=span(\bigcap\limits_{e\in X_1}B_{r(e)})\subseteq span(B_A)$, where the second to last equality follows from Claim 1. Similarly we get $\prod\limits_{f\in X_2}\pi(p_{r(f)})(N)\subseteq span(B_A)$ and $\prod\limits_{e\in X_1}\pi(p_{r(e)})\prod\limits_{f\in X_2}\pi(p_{r(f)})(N)\subseteq span(B_A)$. Therefore $N_A=\pi(p_A)(N)\subseteq span(B_A)$. The general case, that is, the case $A=\left(\bigcap\limits_{e\in X_1}r(e) \right)\cup...\cup\left(\bigcap\limits_{e\in X_n}r(e) \right)\bigcup F$ follows similarly and we leave the details to the reader. 

Since the set $B$ is linearly independent and $B_A\subseteq B$ if follows that $B_A$ is linearly independent an hence $B_A$ is a basis of $N_A$. Furthermore, $B_A\subseteq B$ also implies that $B_A$ is well defined.

The last statement of the lemma follows directly from the definition of the sets $B_A$, $B_D$, $B_{A\cup D}$ and $B_{A\cap D}$, where $A, D \in \GG^0$.
\fim

For the next theorem we recall the following definition:

\begin{definicao}\label{equivrepres} Let $\pi: L_R(\GG) \rightarrow Hom_K(M)$ and $\varphi: L_R(\GG) \rightarrow Hom_R(N)$ be representations of $L_R(\GG)$, where $M$ and $N$ are $R$-modules. We say that $\pi$ is equivalent to $\varphi$ if there exists a $R$-module isomorphism $T:M\rightarrow N$ such that the diagram 
\begin{displaymath}
\xymatrix{
M \ar[r]^{\pi(a)} \ar[d]_{T} &
M \ar[d]^{T} \\
N \ar[r]_{\varphi(a)} & N }
\end{displaymath}
commutes, for each $a\in L_R(\GG)$.
\end{definicao}

It is not true that each representation of $L_R(\GG)$ is equivalent to a representation induced from an $\GG$-algebraic branching system. See, for example, Remark 5.2 of \cite{represleavitt}. However, we get the following theorem:

\begin{teorema}\label{semreferencia} Let $\varphi:L_R(\GG)\rightarrow Hom_R(N)$ be a permutative homomorphism, and let $B$ be a basis of $N$ satisfying the conditions of Definition \ref{permutativerep}. Suppose that $\varphi(p_A)(h_x)=0$ for each $h_x\in B\setminus B_A$ and $A\in \GG^0$, and that $\varphi(s_e^*)(h_x)=0$ for each edge $e$ and $h_x\in B\setminus B_e$. Then there exists a $\GG$-algebraic branching system $X$ such that the representation $\pi:L_R(\GG)\rightarrow Hom_R(M)$, induced by Proposition~\ref{repinducedbybranchingsystems}, is equivalent to $\varphi$, where $M$ is the $R$-module of all the maps from $X$ to $R$ with finite support.\color{red} 
\end{teorema}

\demo Let $B=\{h_x\}_{x\in X}$ be a basis of $N$, with subsets $B_e$, $B_v$ and $B_{r(e)}$ for each edge $e$ and vertex $v$ satisfying the conditions of Definition \ref{permutativerep}. For each $A\in \GG^0$, let $B_A$ be the set defined just before Lemma~\ref{welldefined}. 
For each $e\in \GG ^1$, define $R_e=\{x\in X:h_x\in B_e\}$ and $D_{r(e)}=\{x\in X:h_x\in B_{r(e)}\}$ (notice that $X$ is the index set of $B$). Moreover, for each $A\in \GG^0$ define $D_A=\{x\in X:h_x\in B_A\}$. For a given edge $e$ recall that the map $\varphi(s_e):N_{r(e)}\rightarrow N_e$ is an isomorphism and that $\varphi(s_e)(B_{r(e)})=B_e$. So we get a bijective map $f_e:D_{r(e)}\rightarrow R_e$, such that $\varphi(s_e)(h_x)=h_{f_e(x)}$. 

It is not hard to see that $X$ together with the subsets $\{R_e, D_A\}_{e\in \GG^1, A\in \GG^0}$, and the maps $f_e:D_{r(e)}\rightarrow R_e$ defined above, is a $\GG$-algebraic branching system. For example, to see that $D_{A\cap C}=D_A\cap D_C$ for $A,C\in \GG^0$, notice that $$D_{A\cap C}=\{x\in X:h_x\in B_{A\cap C}\}=\{x\in X:h_x\in B_A\cap B_C\}=D_A\cap D_C,$$ where the second to last equality follows from Lemma \ref{welldefined}, more precisely, from the fact that $B_{A\cap C}=B_A\cap B_C$. Similarly the equality $D_{A\cup C}=D_A\cup D_C$ also holds. The verification of the other conditions of Definition \ref{branchsystem} are left to the reader.

Let $M$ be the $R$ module of all the maps from $X$ to $R$ with finite support, and let $\pi:L_R(\GG)\rightarrow Hom_R(M)$ be the homomorphism induced as in Proposition~\ref{repinducedbybranchingsystems} and Remark \ref{finitesupport}.

Let $\delta_x\in M$ be the map defined by $\delta_x(y)=0$ if $y\neq x$ and $\delta_x(x)=1$. Notice that $\{\delta_x\}_{x\in X}$ is a basis of $M$. Let $T:M\rightarrow N$ be the isomorphism defined by $T(\sum\limits_{i=1}^n k_i\delta_{x_i})=\sum\limits_{i=1}^n k_ih_{x_i}$.

It remains to show that $\varphi(a)=T\circ \pi(a) \circ T^{-1}$, for each $a\in L_R(\GG)$. Notice that for this it is enough to verify that $\varphi(s_e)=T\circ \pi(s_e)\circ T^{-1}$, $\varphi(s_e^*)=T\circ \pi(s_e^*) \circ T^{-1}$, and $\varphi(p_A)=T\circ \pi(p_A) \circ T^{-1}$, for each edge $e$ and $A\in \GG^0$.

Let $A\in \GG^0$. We show that, for each $h_x\in B$, $\varphi(p_A)(h_x)=(T\circ \pi(p_A)\circ T^{-1})(h_x)$. For this, suppose that $h_x\in B_A$. Then $$(T\circ \pi(p_A))(T^{-1}(h_x))=T(\pi(p_A))(\delta_x)=T(\delta_x)=h_x=\varphi(p_A)(h_x).$$ If 
 $h_x\in B\setminus B_A$ then $\varphi(p_A)(h_x)=0$ by hypothesis and, since $x\notin D_A$, we also have $\pi(p_A)(\delta_x)=0$. Hence $(T\circ \pi(p_A))(T^{-1}(h_x))=\varphi(p_A)(h_x)$ for each $h_x\in B$, and therefore $T\circ \pi(p_A)\circ T^{-1}=\varphi(p_A)$.
 
Next we show that $\varphi(s_e)=T\circ\pi(s_e)\circ T^{-1}$. Let $h_x\in B_{r(e)}$. Then $\varphi(s_e)(h_x)=h_{f_e(x)}$, and  $T(\pi(s_e)(T^{-1}(h_x)))=T(\pi(s_e)(\delta_x))=T(\delta_x\circ f_e^{-1})=T(\delta_{f_e(x)})=h_{f_e(x)}$. For $h_x\in B\setminus B_{r(e)}$ we get $\varphi(s_e)(h_x)=\varphi(s_e)\varphi(s_e^*)\varphi(s_e)(h_x)=\varphi(s_e)\varphi(p_{r(e)})(h_x)=0$ (since $\varphi(p_{r(e)})(h_x)=0$ by hypothesis), and $\pi(s_e)(T^{-1}(h_x))=\pi(s_e)(\delta_x)=\delta_x \circ f_e^{-1}.1_{R_e}=0$ since $x\notin D_{r(e)}$.

It remains to prove that $\varphi(s_e^*)=T\circ\pi(s_e^*)\circ T^{-1}$. For $h_x\in B\setminus B_e$ we have that $\varphi(s_e^*)(h_x)=0$ by hypothesis, and $T(\pi(s_e^*)(T^{-1}(h_x)))=T(\pi(s_e^*)(\delta_x))=T(\delta_x\circ f_e.1_{D_{r(e)}})=0$, since $x\notin R_e$. For $h_x\in B_e$ we get $\varphi(s_e^*)(h_x)=h_{f_e^{-1}(x)}$, and $T(\pi(s_e^*)(T^{-1}(h_x)))=T(\pi(s_e^*)(\delta_x))=T(\delta_x\circ f_e.1_{D_{r(e)}})=T(\delta_{f_e^{-1}(x)})=h_{f_e^{-1}(x)}$.
\fim

\begin{remark} 
If $e$ is an edge in $\mathcal{G}$ such that $s(e)$ is a finite emitter then, in the previous theorem, the hypothesis $\pi(s_e^*)(h)=0$ for each $h\in B\setminus B_e$, follows from the hypothesis $\pi(p_A)(h)=0$ for each $h\in B\setminus B_A$. In fact, let $v=s(e)$. Since $p_v=\sum\limits_{s(f)=v}\pi(s_f)\pi(s_f)^*$ (because $s^{-1}(v)$ is finite) then $B_v=\bigcup\limits_{f\in s^{-1}(v)}B_f$, where the last union is a disjoint union. Then, for $h\in B\setminus B_v$, we get $\pi(p_v)(h)=0$ (by hypothesis). For $h\in B_v\setminus B_e$, let $f\in s^{-1}(v)$ be such that $h\in B_f$ and $f\neq e$. Notice that $h=\pi(s_f)\pi(s_f^*)(h)$, and then $\pi(s_e^*)(h)=\pi(s_e^*)\pi(s_f)\pi(s_f^*)(h)=0$, since $\pi(s_f)\pi(s_f^*)=0$.

Therefore, if $\mathcal{G}$ has no infinite emitters, the hypothesis $\pi(s_e^*)(h)=0$ for each $h\in B\setminus B_e$ is unnecessary.
\end{remark}

We now proceed to describe ultragraphs for which a large class of representations is permutative. We recall some definitions and propositions from \cite{ultragraphrep}.

\begin{definicao}[\cite{ultragraphrep}, Definition 6.5]
Let $\mathcal{G}$ be an ultragraph. An \emph{extreme vertex} is an element $A\in r(\mathcal{G}^1)\cup s(\mathcal{G}^1)$ satisfying
\begin{enumerate}
\item either $A=r(e)$ for some edge $e$ and $A\cap r(\mathcal{G}^1\setminus\{e\})=\emptyset=A\cap s(\mathcal{G}^1)$; or
\item $A=s(e)$ for some edge $e$ and $A\cap s(\mathcal{G}^1\setminus \{e\})=\emptyset=A\cap r(\mathcal{G}^1)$.
\end{enumerate}
The edge $e$ associated to an extreme vertex $A$ as above is called the \emph{extreme edge} of $A$.
\end{definicao}

Let $\mathcal{G}$ be an ultragraph. Define the set of \emph{isolated vertices} of $\mathcal{G}$ to be
\[
I_0:=\Big\{v\in G^0: v \notin \Big(\Big(\bigcup_{e \in \mathcal{G}^1}r(e)\Big) \cup s(\mathcal{G}^1)\Big)\Big\}
\]
and define the ultragraph $\mathbb{G}_0:=(G^0 \setminus I_0, \mathcal{G}^1,r,s)$. Denote by $X_1$ the set of extreme vertices of $\mathbb{G}_0$, let $\overline{X_1}=\bigcup\limits_{A\in X_1}A$, and denote by $Y_1$ the set of extreme edges of $\mathbb{G}_0$. Notice that the extreme vertices and the extreme edges of $\mathcal{G}$ and $\mathbb{G}_0$ are the same. Denote by $I_1$ the set of isolated vertices of the ultragraph $\Big(G^0 \setminus (I_0 \cup \overline{X_1}),\mathcal{G}^1\setminus Y_1,r,s\Big)$, and define
$$\mathbb{G}_1=\Big(G^0 \setminus (I_0 \cup I_1\cup \overline{X_1}),\mathcal{G}^1\setminus Y_1,r,s\Big).$$
Now, define $X_2$ and $Y_2$ as being the extreme vertices and extreme edges of the ultragraph $\mathbb{G}_1$, let $\overline{X_2}=\bigcup\limits_{A\in X_2}A$, let $I_2$ be the isolated vertices of the ultragraph $$\Big(G^0 \setminus \big(I_0 \cup I_1\cup \overline{X_1}\cup\overline{X_2}\big),\mathcal{G}^1\setminus (Y_1\cup Y_2),r,s\Big)$$ and let

$$\mathbb{G}_2=\Big(G^0 \setminus \big(I_0 \cup I_1\cup I_2\cup \overline{X_1}\cup\overline{X_2}\big),\mathcal{G}^1\setminus (Y_1\cup Y_2),r,s\Big).$$

Inductively, while $X_n\neq\emptyset$, we define the ultragraphs $\mathbb{G}_n$ and the sets $X_{n+1}$, of extreme vertices of $\mathbb{G}_n$, and $Y_{n+1}$, of extreme edges $\mathbb{G}_n$. We also define the sets $\overline{X_{n+1}}=\bigcup\limits_{A\in X_n}A$ and the set of isolated vertices $I_{n+1}$ of the ultragraph $\mathbb{G}_n$.

Notice that there is a bijective correspondence between the sets $X_n$ and $Y_n$, associating each extreme vertex $A\in X_n$ to an unique extreme edge $e\in Y_n$. For each $A\in X_n$, let $e\in Y_n$ be the (unique) edge associated to $A$. If $A=r(e)$ then $A$ is called a {\it final vertex} of $X_n$ and, if $A=s(e)$, then $A$ is called an {\it initial vertex} of $X_n$. We denote the set of initial vertices of $X_n$ by $X_n^{\mathrm{ini}}$ and the set of final vertices of $X_n$ by $X_n^{\mathrm{fin}}$.

The following theorem is the algebraic version of  \cite[Theorem~6.8]{ultragraphrep}. Until this moment, the coefficient ring $R$ in the Leavitt path algebras of ultragraphs appearing in this paper was assumed only to be a unital commutative ring. However, for general modules over commutative rings, it is not true that each submodule has a basis. But this fact, which we need in the next theorem, is true if $R$ is a field.


\begin{teorema}\label{esteaqui}
Let $\mathcal{G}$ be an ultragraph, $R$ be a field, $N$ be an $R$-module, and let $\pi:L_R(\mathcal{G})\rightarrow Hom_R(N)$ be a representation. Let $N_{r(e)}$ and $N_v$ be as in the beginning of this section. Suppose that $N_{r(e)}=\oplus_{v\in r(e)}N_v$,  for each $e \in \mathcal{G}^1$. If there exists $n \geq 1$ such that $X_1,\dots,X_n \neq \emptyset$, and $\Big(\bigcup\limits_{e \in \mathcal{G}^1}r(e)\Big) \cup s(\mathcal{G}^1)=\bigcup\limits_{i=1}^{n}(\overline{X_i} \cup I_i)$, then $\pi$ is permutative.
\end{teorema}

\demo The proof of this theorem is analogous to the proof of Theorem 6.8 in \cite{ultragraphrep}. The only difference is that in the proof of Theorem 6.8 in \cite{ultragraphrep} the bases are orthonormal bases, while here they are bases only.
\fim

\begin{remark} The ideas of the proof of the previous theorem may be applied to a larger class of ultragraphs than the one satisfying the  hypothesis of the theorem. For example, the ultragraph $\mathcal{G}$

\vspace{2 cm}
\centerline{
\setlength{\unitlength}{1.5cm}
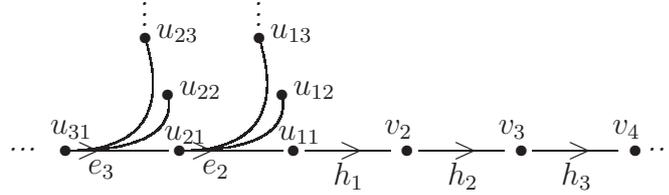
\begin{picture}(14,0)
\put(10,0){\circle*{0.1}}
\put(10.1,0){...}
\put(9.8,0.15){$v_4$}
\put(9,0){\circle*{0.1}}
\put(8.8,0.15){$v_3$}
\put(8,0){\circle*{0.1}}
\put(7.8,0.15){$v_2$}
\put(9.1,0){$\line(1,0){0.75}$}
\put(9.4,-0.07){$>$}
\put(8.1,0){$\line(1,0){0.75}$}
\put(8.4,-0.07){$>$}
\put(7,0){\circle*{0.1}}
\put(7.1,0){$\line(1,0){0.75}$}
\put(7.4,-0.07){$>$}
\put(6.87,0.1){$u_{11}$}
\put(6,0){\circle*{0.1}}
\qbezier(6,0)(7,0)(6.9,0.5)
\qbezier(6,0)(6,0)(6.9,0)
\qbezier(6,0)(7,0)(6.7,1)
\put(6.9,0.5){\circle*{0.1}}
\put(7,0.5){$u_{12}$}
\put(6.7,1){\circle*{0.1}}
\put(6.8,1){$u_{13}$}
\put(6.1,-0.07){$>$}
\put(6.65,1.1){\vdots}
\put(5,0){\circle*{0.1}}
\qbezier(5,0)(6,0)(5.9,0.5)
\qbezier(5,0)(5,0)(5.9,0)
\qbezier(5,0)(6,0)(5.7,1)
\put(5.9,0.5){\circle*{0.1}}
\put(6,0.5){$u_{22}$}
\put(5.7,1){\circle*{0.1}}
\put(5.8,1){$u_{23}$}
\put(5.1,-0.07){$>$}
\put(5.65,1.1){\vdots}
\put(5.87,0.1){$u_{21}$}
\put(4.87,0.15){$u_{31}$}
\put(4.5,0){...}
\put(5.2,-0.2){$e_3$}
\put(6.2,-0.2){$e_2$}
\put(7.35,-0.3){$h_1$}
\put(8.35,-0.3){$h_2$}
\put(9.35,-0.3){$h_3$}
\end{picture}}
\vspace{1cm}
 does not satisfy the hypothesis of the previous theorem, but each representation $\pi:L_R(\mathcal{G})\rightarrow Hom_R(N)$ is permutative. See [Remark 6.9, \cite{ultragraphrep}] for more details.
\end{remark}

We end the paper with the following result, regarding unitary equivalence of representations.

\begin{teorema}
Let $\mathcal{G}$ be an ultragraph and suppose that there exists $n \geq 1$ such that $X_1,\dots,X_n \neq \emptyset$, and $\Big(\bigcup\limits_{e \in \mathcal{G}^1}r(e)\Big) \cup s(\mathcal{G}^1)=\bigcup\limits_{i=1}^{n}(\overline{X_i} \cup I_i)$. Let $R$ be a field, $N$ be an $R$-module, and let $\pi:L_R(\mathcal{G})\rightarrow Hom_R(N)$ be a representation. Let $N_{r(e)}$, $N_e$ and $N_v$ be as in the beginning of this section. Suppose that $N_{r(e)}=\oplus_{v\in r(e)}N_v$, for each $e \in \mathcal{G}^1$. Let $M=\pi(L_R(\mathcal{G}))(N)$, and let $\widetilde{\pi}:L_R(\mathcal{G})\rightarrow Hom_R(M)$ be the restriction of $\pi$. If $N_v=\bigoplus_{e\in s^{-1}(v)}N_e$, for each vertex $v$ which is not a sink, then $\widetilde{\pi}$ is equivalent to a representation induced by a branching system.
\end{teorema}

\demo By Theorem \ref{esteaqui} we get that $\widetilde{\pi}$ is permutative. Let $B$ be the basis obtained in the proof of Theorem~\ref{esteaqui}, that is, $B$ is a basis of $M$, $B_v$ is a basis of $\widetilde{\pi}(p_v)(M)$ for each vertex $v$, $B_e$ is a basis of $\widetilde{\pi}(s_e)\widetilde{\pi}(s_e^*)(M)$ for each edge $e$, $B_v\supseteq B_e$ for each $e\in s^{-1}(v)$, and $B_v\subseteq B$ for each vertex $v$. Moreover, by hypothesis, $B_{r(e)}=\bigcup\limits_{v\in r(e)}B_v$. Notice that $B=\bigcup\limits_{v\in G^0}B_v$.   
By Theorem~\ref{semreferencia}, we need to show that $\widetilde{\pi}(p_A)(h_x)=0$ for each $h_x\in B\setminus B_A$, and $\widetilde{\pi}(s_e^*)(h_x)=0$ for all $h_x\in B\setminus B_e$. 

Let $h_x\in B\setminus B_A$. Then $h_x\in B_u$ for some $u\notin A$, and so $\widetilde{\pi}(p_A)(h_x)=\widetilde{\pi}(p_A)\widetilde{\pi}(p_u)(h_x)=0$ ( $\widetilde{\pi}(p_A)\widetilde{\pi}(p_u)=0$).

Let $e$ be an edge and $h_x\in B\setminus B_e$. If $h_x\in B_u$ with $u\neq s(e)$ then $\widetilde{\pi}(s_e^*)(h_x)=\widetilde{\pi}(s_e^*)\widetilde{\pi}(p_{s(e)})(h_x)=0$ since $\widetilde{\pi}(p_{s(e)})(h_x)=0$. So, let $h_x\in B_u$ with $u=s^{-1}(e)$. Since $B_u=\bigcup\limits_{f\in s^{-1}(u)}B_f$, there exists $f\in s^{-1}(u)$ with $f\neq e$ such that $h_x\in B_f$. From the proof of the previous theorem we get that $B_f=\widetilde{\pi}(s_f)(B_{r(f)})$, and so $h_x=\widetilde{\pi}(s_f)(h)$ for some $h\in B_{r(f)}$. Then $\widetilde{\pi}(s_e^*)(h_x)=\widetilde{\pi}(s_e^*)\widetilde{\pi}(s_f)(h)=0$ since $\widetilde{\pi}(s_e^*)\widetilde{\pi}(s_f)=0$.

\fim

\begin{remark}  Note that in the previous theorem, the condition $N_v=\bigoplus_{e\in s^{-1}(v)}N_e$ is automatically satisfied if $0<|s^{-1}(v)|<\infty$.
\end{remark}

\vspace{1.5pc}

Daniel Gon\c{c}alves, Departamento de Matem\'{a}tica, Universidade Federal de Santa Catarina, Florian\'{o}polis, 88040-900, Brasil

Email: daemig@gmail.com

\vspace{0.5pc}
Danilo Royer, Departamento de Matem\'{a}tica, Universidade Federal de Santa Catarina, Florian\'{o}polis, 88040-900, Brasil

Email: danilo.royer@ufsc.br
\vspace{0.5pc}

\end{document}